\documentclass[11pt]{amsart}
\usepackage{amsmath,amsfonts,amscd,amssymb}

\theoremstyle{plain}

\newtheorem{theorem}[equation]{Theorem}
\newtheorem{conjecture}[equation]{Conjecture}
\newtheorem{proposition}[equation]{Proposition}
\newtheorem{lemma}[equation]{Lemma}
\newtheorem{corollary}[equation]{Corollary}
\theoremstyle{definition}

\newtheorem{remark}[equation]{Remark}

\def\O{{\mathcal O}}
\def\X{{\mathcal X}}

\def\F{{\mathbb F}}

\def\Q{{\mathbb Q}}

\def\R{{\mathbb R}}

\def\Z{{\mathbb Z}}

\def\triv{{\mathbf 1}}

\def\rksel#1#2#3{\rk_{#3}(#1/#2)}
\def\rkan#1#2{\ord_{s=1}L(#1/#2,s)}

\def\newmathop#1{\expandafter\gdef\csname #1\endcsname{\mathop{\rm #1}\nolimits}}
\newmathop{Tr}

\newmathop{rk}
\newmathop{ord}
\newmathop{Gal}
\newmathop{BSD}
\newmathop{Hom}
\newmathop{tors}
\newmathop{coker}
\newmathop{Ind}
\newmathop{Res}
\newmathop{Frob}
\let\oldchar\char\newmathop{char}\let\char\oldchar
\newmathop{GL}
\newmathop{PGL}
\newmathop{SL}
\newmathop{Aut}
\newmathop{Sel}
\newmathop{End}
\newmathop{Maps}
\newmathop{Re}
\newmathop{Reg}
\newmathop{Sy}
\newmathop{Syl}
\newmathop{sgn}

\let\iso\cong
\let\tensor\otimes

\def\beq{$$\begin{array}{llllllllllllllll}}
\def\eeq{\end{array}$$}
\def\beqn{\begin{equation}\begin{array}{llllllllllllllll}}
\def\eeqn{\end{array}\end{equation}}

\input cyracc.def       
\font\tencyr=wncyr10
\def\sha{\text{\tencyr\cyracc{Sh}}}

\def\qequal{\hbox{\rlap{$\,\raise-1pt\hbox{?}$}=}}

\def\tbuildrel#1\over#2{\buildrel\text{\rm\normalsize\smaller[3]#1}\over{#2}}

\def\thincdots{\raise1.3pt\hbox{$\scriptscriptstyle
  \>\cdot\>\cdot\>\cdot\>\cdot\hskip0.3pt$}}

\def\neron#1{\omega_{#1}^o}

\def\subgroup{\raise0.5pt\hbox{$\,\scriptstyle<\,$}}
\let\<\subgroup

\theoremstyle{plain}
\newtheorem*{conjectures}{Conjectures}

\def\rksel#1#2#3{\rk_{#3}#1/#2}
\let\p\wp
\def\m{{\mathfrak m}}
\begin{document}

\let\introdagger\dagger
\title{A note on Larsen's conjecture and ranks of elliptic curves}
\author{Tim$^\introdagger$ and Vladimir Dokchitser}
\date{March 6, 2008}
\thanks{$^\dagger$Supported by a Royal Society University Research Fellowship}

\address{Robinson College, Cambridge CB3 9AN, United Kingdom}
\email{t.dokchitser@dpmms.cam.ac.uk}
\address{Gonville \& Caius College, Cambridge CB2 1TA, United Kingdom}
\email{v.dokchitser@dpmms.cam.ac.uk}
\maketitle

\section{Introduction}

In \cite{Lar} M. Larsen proposed the following:

\begin{conjecture}
\label{larsen}
Let $K$ be a finitely generated infinite field, $K_s$ its separable closure,
$G$ a finitely generated subgroup of $\Gal(K_s/K)$, and $K_{\!s}{}\!^G$ its
fixed field.
If $E/K$ is an elliptic curve, then the group
$E(K_{\!s}{}\!^G)\tensor\Q$ is infinitely generated.
\end{conjecture}

When $G$ is generated by one element, this conjecture is known
for $K=\Q$ \cite{Im} and in a number of other cases \cite{BI,Im2,Im3};
see also \cite{Loz} for some results when $G$ has two generators.
We give evidence for the full conjecture for $K=\Q$
by showing that it is
implied both by the Shafarevich-Tate
and the (first part of) the Birch--Swinnerton-Dyer conjecture
for elliptic curves over number fields:

\begin{theorem}
\label{main1}
Let $E/\Q$ be an elliptic curve and $G\subset \Gal(\bar\Q/\Q)$ a finitely
generated subgroup.
The analytic rank and the $p^\infty$-Selmer rank for every odd $p$ are
unbounded in number fields contained in $\bar{\Q}^G$.
\end{theorem}

Recall that for an elliptic curve $E$ over a number field $F$,
the {\em analytic rank\/} $\rk_{an} E/F$ is defined as $\ord_{s=1}L(E/F,s)$.
Even for $E$ defined over $\Q$
the analytic continuation of $L(E/F,s)$ to $s=1$ is not known for general $F$,
but in the proof of the theorem we will choose number fields where it is.
The {\em $p^\infty$-Selmer rank\/} $\rk_p E/F$ is the usual Mordell-Weil
rank $\rk E/F=\dim E(F)\tensor\Q$ plus the $\Z_p$-corank of $\sha(E/F)$,
the number of copies of $\Q_p/\Z_p$ in $\sha(E/F)$.
We remind the reader of the usual conjectures concerning these ranks:


\begin{conjectures}
Let $K$ be a number field, and $E/K$ an elliptic curve.
\begin{itemize}
\labelsep 0.6cm
\item[$\bullet$]\hspace{-6mm}
$L(E/K,s)$ is entire, and $\rk E/K\!=\!\rk_{an}\!E/K$ (Birch--Swinnerton-Dyer).
\item[$\bullet$]\hspace{-6mm}
$\sha(E/K)$ is finite; thus $\rk E/K\!=\!\rksel EKp$ for all $p$ (Shafarevich-Tate).
\item[$\bullet$]\hspace{-6mm}
$(-1)^{\rk E/K}=w(E/K)$ (Parity conjecture).
\end{itemize}
\end{conjectures}

\noindent
The global root number $w(E/K)$ is defined as a product of local root numbers
\pagebreak\noindent
$w(E/K_v)=\pm 1$ over all places of $K$, and is the conjectural sign in the
functional equation for $L(E/K,s)$. Thus it determines $\rk_{an}(E/K)\mod 2$,
so the parity conjecture is a weaker version of the Birch--Swinnerton-Dyer
conjecture.

In fact,
the proof of Theorem \ref{main1} shows that the parity conjecture
also implies Conjecture \ref{larsen} for elliptic curves over $K=\Q$.
When $K$ is a general number field, the analytic continuation of $L(E/K,s)$
to $s=1$ is not known so we cannot say anything about $\rk_{an}E/K$.
However, we still have the following

\begin{theorem}
\label{main2}
Let $K$ be a number field, $E/K$ an elliptic curve and
$G$ a finitely generated subgroup of $\Gal(\bar K/K)$.
Then both the parity and the Shafarevich-Tate conjectures imply
Conjecture \ref{larsen}, provided
\begin{itemize}
\item $K$ has a real place, or
\item $E$ has non-integral $j$-invariant.
\end{itemize}
\end{theorem}

\begin{remark}
Slightly more generally, the theorem applies if
there is a place $v$ of $K$ and some quadratic extension $L_w/K_v$
where $w(E/L_w)=-1$.
(See \cite{Evilquad} for a classification of such $E/K_v$.)
\end{remark}

The proofs of both theorems rely on an elementary root number argument
in $\F_p^r\rtimes C_2$-extensions
and some version of the parity conjecture.
For the second theorem we will need
the following result, which is a slight variation of \cite{Squarity} Thm. 1.3.

\begin{proposition}
\label{gensq}
Let $K$ be a number field, $E/K$ an elliptic curve and $M/K$
a quadratic extension.
Suppose for every prime $w|6$ of $M$ which is not split in $M/K$,
$E$ has semistable reduction at $w$, and the reduction is not good
supersingular if $w|2$.
If $\sha(E/M(E[2]))[6^\infty]$ is finite, then $$(-1)^{\rk E/M}=w(E/M).$$
\end{proposition}

This has a consequence, which is of interest in itself:

\begin{corollary}
\label{useful}
Let $M/K$ be a quadratic extension of number fields, $E/K$ an elliptic curve
and assume that $\sha(E/M(E[2]))$ is finite.
If all primes of bad reduction of $E$ split in $M/K$,
then the parity conjecture holds for $E/M$,
$$
  (-1)^{\rk E/M} = w(E/M) 
    = (-1)^{\text{{\rm\#Archimedean places of $M$}}}.
$$
\end{corollary}


%
%
%
%
%
%
%

\medskip
\noindent
{\bf Notation.}
Throughout the note $K$ denotes a number field and $E$ an elliptic curve
defined over $K$. We write $w(E/K_v)=\pm 1$ for the local root number of $E$
at $v$, and $w(E/K)=\prod_{v} w(E/K_v)$ for the global root number.
For an Artin representation $\tau$ of $\Gal(\bar K/K)$ we write
$w(E/K_v,\tau)$ and $w(E/K,\tau)=\prod_{v} w(E/K_v,\tau)$ for the local
and global root numbers of the twist of $E$ by~$\tau$.
See e.g. \cite{RohG} 
for properties of root numbers of elliptic curves and their twists.

\bigskip
\noindent
{\bf Acknowledgments.} We would like to thank A. Jensen 
for helpful discussions.

\newpage

\section{Proof of Theorem \ref{main1}}

\begin{lemma}
\label{lem}
Let $M/K$ be a quadratic extension of number fields,
$E/K$ an elliptic curve and
$G\subset \Gal(\bar K/K)$ a finitely generated subgroup.
For every odd prime $p$ and $r\ge 1$ there exists a Galois extension
$F/K$ containing $M$ such that
\begin{enumerate}
\item
$\Gal(F/K)\iso\F_p^r\rtimes C_2$, with $C_2$ acting by~$-1$.
\item
The image of $G$ in $\Gal(F/K)$ has order at most 2.
\item
The primes of $M$ above primes of bad reduction
for $E/K$ split completely in $F$.
\end{enumerate}
For such an extension, $w(E/K,\rho)=w(E/M)$ for every irreducible
2-dimen\-sional representation $\rho$ of $\Gal(F/K)$.
\end{lemma}

\begin{proof}
Pick primes $\p_1,\ldots,\p_n$ of $K$ that split completely in
$M(\zeta_p)$, and consider $\m=\prod_i\p_i$ as a modulus of $M$
($\zeta_p$ denotes a primitive $p$th root of unity).
Write $I^\m$ for the group of fractional ideals of $M$
that are coprime to~$\m$, and $P^\m$ for the subgroup of principal ideals
that can be generated by an element congruent to 1 mod $\m$.

Let $Q$ be the largest $\F_p$-vector space quotient of the group $I^\m/P_\m$.
Then
$2n-\delta \le \dim Q \le 2n$ with
$\delta=\dim\O_M^*/\O_M^{*p}+\dim{\rm Cl}_M[p]$.
Comparing this with the corresponding group of $K$, we see that
the $\Gal(M/K)$-antiinvariant part of $Q$ has dimension $d\ge n-\delta$.
By class field theory it yields a Galois extension $F_n/K$ with
$\Gal(F_n/K)\iso\F_p^d\rtimes C_2$, $C_2$ acting by~$-1$.

Consider $H=\Gal(F_n/M)=\F_p^d$. The group $G\cap\Gal(\bar M/M)$ is of index
at most 2 in $G$, and the $\F_p$-dimension of its image in $H$ is bounded by the
number of generators of $G$. Also, the decomposition subgroup in $H$
of any prime of $M$ has size bounded by a constant independent of $n$.
Now let $H_0\<H$ be generated by the image of $G\cap\Gal(\bar M/M)$ and
the decomposition subgroups of primes of $M$ that lie above primes of bad
reduction for $E/K$. Then $F_n^{H_0}$ satisfies (2) and (3), and letting
$n\to \infty$ gives the desired extensions $F/K$.

For the last claim, let $\epsilon$ be the non-trivial character of
$\Gal(M/K)$. The complex irreducible representations of $\Gal(F/K)$ are
$\triv$, $\epsilon$ and 2-dimensional representations each of which factors
through a $D_{2p}$-quotient.

Recall that $w(E/M)=w(E/K,\triv\oplus\epsilon)$ by inductivity
of global root numbers, so it suffices to check that
$w(E/K_v,\triv\oplus\epsilon)=w(E/K_v,\rho)$ for all places $v$ of $K$.
This holds for places $v$ whose decomposition group $D$ in $F/K$ has order
at most 2, because $\rho$ and $\triv\oplus\epsilon$ are isomorphic
as $D$-representations. In particular, this includes Archimedean places and
places of bad reduction for $E$.
On the other hand, if $E$ has good reduction at $v$, then
for every 2-dimensional self-dual Artin representation $\tau$
of $\Gal(\bar K/K)$,
$$
  w(E/K_v,\tau) = w(\tau/K_v)^2 = \det(\tau(-1)),
$$
by the unramified twist formula \cite{TatN}~3.4.6 and the determinant
formula \cite{TatN}~3.4.7. Here we implicitly use the local reciprocity map
at $v$ to evaluate $\tau$ at~$-1$. As
$\det\rho=\epsilon=\det(\triv\oplus\epsilon)$,
the claim follows.
\end{proof}

We now prove Theorem \ref{main1}.
Take an imaginary quadratic field $M\!=\!\Q(\sqrt{\!-d})$ where all bad primes
for $E$ split. Note that $w(E/M)=-1$, as the contribution from bad primes
is $(\pm 1)^2$, and the local root number is $+1$ for primes of good reduction
and $-1$ for infinite places.

Let $p$ be an odd prime.
By the $p$-parity conjecture for $E/\Q$ and its quadratic
twist by~$-d$ (\cite{Squarity} Thm. 1.4),
we have $(-1)^{\rk_p(E/M)}=w(E/M)$, so $\rk_p(E/M)$ is odd. Let $F$ be as in Lemma \ref{lem} for some $r\ge 1$,
and let $V\subset\F_p^r$ be any index $p$ subgroup, so
$V\triangleleft\Gal(F/\Q)$ and $\Gal(F^V/\Q)\iso D_{2p}$.
The image of $G$ in this Galois group has order at most 2, so there is
a degree $p$ extension $L/\Q$ inside $F^V$ fixed by $G$.
As all bad primes for $E$ split in $M$, by \cite{Squarity} Prop. 4.17,
$$
  \rksel E{M}p+\tfrac{2}{p-1}(\rksel E{L}p-\rksel E{\Q}p)
  \equiv
  0
  \pmod2.
$$
In particular, $\rksel E{L}p>\rksel E{\Q}p$.

Write $X$ for the dual $p^\infty$-Selmer group
$\Hom(\Sel_{p^\infty}(E/F),\Q_p/\Z_p)$ and $\X=X\tensor_{\Z_p}\Q_p$.
Thus $\X$ is a $\Q_p$-valued representation of $\Gal(F/\Q)$, and
$\rksel Ekp=\dim\X^{\Gal(F/k)}$ for every $k\subset F$,
see e.g. \cite{Squarity} Lemma 4.14.

The irreducible $\Q_p$-representations of $\Gal(F^V/\Q)\iso D_{2p}$
are trivial $\triv$, sign $\epsilon$ and $(p-1)$-dimensional irreducible
$\rho_V$. Their invariants under an element of order 2 of $D_{2p}$ are
$1$-, $0$- and $\tfrac{p-1}2$-dimensional, respectively. It follows that
$\X$ contains a copy of $\rho_V$ for every $V$, and $\dim \X^G$ is therefore
at least $\tfrac{p-1}2$ times the number of hyperplanes $V\subset\F_p^r$ (which is
$\tfrac{p^r-1}{p-1}$). Letting $r\to\infty$ we deduce that the
$p^\infty$-Selmer rank of
$E$ is unbounded for number fields inside $\bar\Q^G$.

To prove that the analytic rank is unbounded, let $p=3$ for simplicity and
take the same $F$ as above.
The irreducible complex representations of $\Gal(F/\Q)$ are $\triv, \epsilon$,
and $\rho_V$ for varying $V$. The curve $E/\Q$ is modular, so the
$L$-functions $L(E/\Q,s)$, $L(E,\epsilon,s)$ and $L(E,\rho_V,s)$ are analytic
and satisfy the expected functional equation. (The twists by $\rho_V$ are
Rankin-Selberg products.) 
The same applies to $L(E/k,s)$ for every $k\subset F$; indeed,
by Artin formalism,
$$
  L(E/k,s) \>=\> L(E,\Ind_k^\Q\triv,s) \>=\>
    \!\!\!\!\prod_{\scriptscriptstyle {\tau\in\{\triv,\epsilon\}\cup\{\rho_V\}_V}}\!\!\!\!
    L(E,\tau,s)^{\langle \tau, \Ind_k^\Q\triv \rangle},
$$
with $\Ind_k^\Q$ a shorthand for $\Ind_{\Gal(F/k)}^{\Gal(F/\Q)}$.
Finally, $w(E,\rho_V)=-1$ by Lemma~\ref{lem}, so each $L(E,\rho_V,s)$
vanishes at $s=1$. Hence
$$
  \rkan Ek\ge\frac{[k:\Q]-2}2
$$
as $\Ind_k^\Q\triv$ contains at most one copy of $\triv$ and $\epsilon$.
Now take $k=F^G$ and let $r\to\infty$ as before.

\section{Proof of Theorem \ref{main2}}

We claim there is a quadratic extension $M/K$ where
$w(E/M)=-1$ and $\rk E/M$ is odd. This will suffice, as we can
then use exactly the same const\-ruction as in the proof of
Theorem~\ref{main1} with $K$ in place of $\Q$ and any odd~$p$.

Fix a place $v$ of $K$ which is either real or with $\ord_v j(E)<0$.
We take $M/K$ to be any quadratic extension where
\begin{itemize}
\item $v$ becomes complex if $K_v=\R$.
\item There is a unique prime $v'$ above $v$, and $E$ has split
multiplicative reduction at $v'$ if $\ord_v j(E)<0$.
\item If $w\ne v$ is either Archimedean, a place of bad reduction of $E$ or
divides 2 then $w$ splits in $M/K$.
\end{itemize}
Such a quadratic $M$ exists by the weak approximation theorem,
as we prescribe only its local behaviour at finitely many places.
The existence of an appropriate local extension in the second case follows
from the theory of the Tate curve, see e.g. \cite{Sil2} Thm 5.3.

Write $v'$ for the unique prime of $M$ above $v$. Then $w(E/M_{v'})=-1$
(see e.g. \cite{RohG} Thm. 2) and all other Archimedean places and bad
primes for $E$ contribute $(\pm 1)^2=1$, so $w(E/M)=-1$.
Finally $\rk E/M$ is odd, by either the assumed parity conjecture or
the Shafarevich-Tate conjecture together with Proposition \ref{gensq}.

%
%

\section{Proof of Proposition \ref{gensq}}

\begin{lemma}
\label{lemma2}
Let $K$ be a number field and $E/K$ an elliptic curve with
a $K$-rational 2-torsion point.
Let $\Sigma$ be the set of primes $v|2$ of $K$ where $E$ has additive
or good supersingular reduction. Suppose
$M/K$ is a quadratic
extension where every prime $v\in\Sigma$ either splits or becomes
a place of multiplicative or good ordinary reduction for $E$.
Then $(-1)^{\rksel EM2}=w(E/M)$.
\end{lemma}

\begin{proof}
Choose a model for $E/K$ of the form
$$
  E: \>\> y^2 = x^3 + ax^2 + bx \>,
$$
and let $\phi:E\to E'$ be the $2$-isogeny with $(0,0)$ in the kernel.

Let $\sigma_w(E/M)=(-1)^{\ord_2\coker\phi_w-\ord_2\ker\phi_w}$, where
$\phi_w: E(M_w)\to E'(M_w)$ is the induced map on local points.
By Cassels' formula (see \cite{Isogroot} \S1),
$$
  (-1)^{\rksel EM2}=\prod_w \sigma_w(E/M).
$$
Write $(\cdot,\cdot)_w=\pm 1$ for the Hilbert symbol at $w$. Then
$$
  \prod_{w|v} w(E/M_w) = \prod_{w|v} \sigma_w(E/M) (a,-b)_w(-2a,a^2-4b)_w
$$
for every place $v$ of $K$. Indeed,
this trivially holds for all $v$
that split in $M/K$, and it holds at all other primes by \cite{Isogroot},
proof of Thm. 4 in \S7. The result follows by the product formula
for Hilbert symbols.
\end{proof}

\begin{lemma}
\label{lemma3}
Let $K$ be a number field and $E/K$ an elliptic curve.
Let $\Sigma$ be the set of primes $v|6$ of $K$ where $E$ has additive
reduction.
Suppose $M/K$ is a quadratic
extension where every prime $v\in\Sigma$ either splits or becomes
a place of semistable reduction for $E$.
Let $F/K$ be any $S_3$-extension not containing $M$, and
write $N=F^{C_2}M$, $N'=F^{C_3}M$. Then
$$
  w(E/M)w(E/N)w(E/N') = (-1)^{\rksel EM3+\rksel EN3+\rksel E{N'}3}.
$$
\end{lemma}

\begin{proof}
Fix a global differential $\omega\ne 0$ for $E/K$. For an extension
$k/K$ and a prime $w$ of $k$ write
$C_w(E/k) = c_w(E/k) |{\omega}/{\neron{w}}|_w$, where
$c_w$ is the local Tamagawa number, $\neron{w}$ a N\'eron differential at $w$,
and $|\cdot|_w$ the normalised $w$-adic absolute value.
By \cite{Squarity} Thm. 4.11 with $p=3$,
$$
  \rksel EM3+\rksel EN3+\rksel E{N'}3 \equiv
  \ord_3
    \frac{\prod_z C_z(E/MF)}{\prod_u C_u(E/N')}
  \mod 2,
$$
the products taken over the primes of $MF$ and $N'$ respectively.
It suffices to show that for every place $v$ of $K$,
$$
  \prod\limits_{s|v} w(E/M_s) \prod\limits_{t|v} w(E/N_t) \prod\limits_{u|v} w(E/N'_u)
    =
  \frac{\prod_{z|v}(-1)^{\ord_3 C_z(E/MF)}}{\prod_{u|v}(-1)^{\ord_3 C_u(E/N')}}\>,
$$
where we interpret the right-hand side as 1 for Archimedean places.
If $v$ splits in $M/K$ the formula trivially holds as each term occurs an
even number of times. For all other $v$, this is proved in
\cite{Squarity} Prop 3.3
(the $G$-set argument in Case 1 also covers Archimedean places).
\end{proof}

We now prove Proposition \ref{gensq}.

The assumption on $\sha$ also forces $\sha(E/k)[6^\infty]$ to be finite
for all intermediate fields $K\subset k\subset M(E[2])$,
see e.g. \cite{Squarity} Rmk. 2.10.
Now we proceed as in the proof of \cite{Squarity} Thm. 3.6.

Write $F=K(E[2])$.
If $E(K)[2]\ne 0$, apply Lemma \ref{lemma2}.
Otherwise $\Gal(F/K)$ is either $C_3$ or $S_3$. In the former case,
$\Gal(FM/M)\iso C_3$ as well, thus $\rk E/M$ and
$\rk E/FM$ have the same parity. It is also well-known that
global root numbers are unchanged in odd degree cyclic extensions
(this follows from \cite{TatN} 3.4.7, 4.2.4),
so $w(E/FM)=w(E/M)$ and the result again follows from Lemma \ref{lemma2}
applied to $FM/F$.

In the last case $\Gal(F/K)\iso S_3$, write $N=F^{C_2}M$, $N'=F^{C_3}M$.
If $M\not\subset F$, the above argument shows that $w(E/N)=(-1)^{\rk E/N}$
and similarly for $N'$; now apply Lemma \ref{lemma3}.
If $M\subset F$, then $F/M$ is a Galois cubic extension, so
we may again show that $w(E/F)=(-1)^{\rk E/F}$.
But this holds by Lemma \ref{lemma2} applied to the quadratic extension $F/N$.

\smallskip

Finally, Corollary \ref{useful} is immediate from the fact that $w(E/K_v)=-1$
for Archimedean $v$.

\end{document}